\documentclass{amsart}

\usepackage{amsmath,amssymb,latexsym}

\numberwithin{equation}{section}
\newtheorem{theorem}{Theorem}[section]

\newtheorem{definition}[theorem]{Definition}

\newtheorem{lemma}[theorem]{Lemma}

\newtheorem{example}[theorem]{Example}

\begin{document}

\pagenumbering{arabic}
\pagestyle{headings}
\def\sof{\hfill\rule{2mm}{2mm}}
\def\SS{\mathcal S}
\def\qq{{\bold q}}
\def\txx{{\frac1{2\sqrt{x}}}}
\def\mn{\text{-}}

\title{{\sc Words restricted by $3$-letter generalized multipermutation patterns}}

\author{Alexander Burstein \and Toufik Mansour}

\date{December 27, 2001}

\address{Department of Mathematics, Iowa State University, Ames, IA
50014-2064 USA}

\email{burstein@math.iastate.edu}

\address{LaBRI, Universit\'e Bordeaux, 351 cours de la
Lib\'eration, 33405 Talence Cedex, France}

\email{toufik@labri.fr}

\maketitle

%===========================================================================

\begin{abstract}
We find exact formulas and/or generating functions for the number
of words avoiding 3-letter generalized multipermutation patterns
and find which of them are equally avoided.
\end{abstract}

%===========================================================================
\section{Introduction}

A \emph{generalized pattern} $\tau$ is a (possibly hyphenated)
string in $[\ell]^m$ which contains all letters in
$[\ell]=\{1,\dots,\ell\}$. We say that the string $\sigma\in[k]^n$
\emph{contains} a generalized pattern $\tau$ if $\sigma$ contains
a subsequence isomorphic to $\tau$ in which the entries
corresponding to consecutive entries of $\tau$ not separated by a
hyphen must be adjacent. Otherwise, we say that $\sigma$
\emph{avoids} $\tau$ and write $\sigma\in [k]^n(\tau)$. Thus,
$[k]^n(\tau)$ denotes the set of strings in $[k]^n$ (i.e. $n$-long
$k$-ary strings) which avoid $\tau$.
\begin{example}
\rm An string $\pi=a_1a_2\dots a_n$ avoids $13\mn2$ if $\pi$ has
no subsequence $a_ia_{i+1}a_j$ with $j>i+1$ and $a_i<a_j<a_{i+1}$.
\end{example}

Classical patterns are generalized patterns with all possible
hyphens (say, $2\mn1\mn3$), in other words, those that place no
adjacency requirements on $\sigma$. The first case of classical
patterns studied was that of permutations avoiding a permutation
pattern of length 3. Knuth \cite{Knuth} found that, for any
$\tau\in S_3$, $|S_n(\tau)|=C_n$, the $n$th Catalan number. Later,
Simion and Schmidt \cite{SS} determined the number $|S_n(P)|$ of
permutations in $S_n$ simultaneously avoiding any given set of
patterns $P\subseteq S_3$. Burstein \cite{Burstein} extended this
to $|[k]^n(P)|$ with $P\subseteq S_3$. Burstein and Mansour
\cite{BM} considered forbidden patterns with repeated letters.

Generalized permutation patterns were introduced by Babson and
Steingr\'{\i}msson \cite{BS} with the purpose of the study of
Mahonian statistics. Later, Claesson \cite{Claesson} and Claesson
and Mansour \cite{CM} considered the number of permutations
avoiding one or two generalized patterns with one hyphen. Finally,
in this paper, we consider the case of words (repeated letters
allowed) avoiding a single generalized pattern of length 3 (with
repeated letters allowed as well).

\begin{definition} \label{defWilf}
We say that two patterns $\tau_1$ and $\tau_2$ are
\emph{Wilf-equivalent} or belong to the same \emph{Wilf class} if
$|[k]^n(\tau_1)|=|[k]^n(\tau_2)|$ for all integers $k,n\ge 0$.
\end{definition}

Given a generalized pattern $\tau$ we define its \emph{reversal}
$r(\tau)$ to be $\tau$ read right-to-left, including hyphens (e.g.
example, $r(13\mn2)=2\mn31$). We also define the \emph{complement}
of $\tau$, denoted $c(\tau)$, to be the pattern obtained by
substituting $\ell+1-\tau(i)$ for $\tau(i)$ and leaving hyphens in
the same positions (e.g., $c(13\mn2)=31\mn2$). Clearly, $c\circ
r=r\circ c$, so $\langle c,r \rangle=D_2$, the group of symmetries
of a rectangle. We call the set $\{\tau, r(\tau), c(\tau),
c(r(\tau))\}$ the \emph{symmetry class} of $\tau$. Obviously, all
patterns in the same symmetry class (e.g. 13-2, 2-31, 31-2, 2-13)
are Wilf-equivalent, so we only need to consider one
representative of each symmetry class to determine Wilf classes.

Let $\tau$ be any generalized pattern; we define $\#\tau(\sigma)$
as the number of occurrences of $\tau$ in $\sigma$. Let
$F_\tau(n,k;q;a_1,\dots,a_d)$ be the {\em occurrence polynomial},
whose coefficient at $q^r$ is the number of words in $[k]^n$
having exactly $r$ occurrences of $\tau$ and ending on $a_1\dots
a_d$; that is,
\[
F_\tau(n,k;q;a_1,\dots,a_d)=\sum_{\sigma\in [k]^{n-d}}
q^{\#\tau(\sigma a_1\dots a_d)}.
\]
For $d=0$, we denote $F_\tau(n,k;q)=F_\tau(n,k;q;\emptyset)$. We
also denote the generating function for the sequence
$\{F_\tau(n,k;q)\}_{n\ge 0}$ by $F_\tau(x;k;q)$; that is,
\[
F_\tau(n,k;q)=\sum_{n\ge 0}x^nF_\tau(n,k;q)=\sum_{n\ge
0}\sum_{\sigma\in [k]^n} q^{\#\tau(\sigma)}x^n.
\]

%===========================================================================
\section{Two-letter generalized patterns}

Here, the only symmetry classes with repeated letters are those of
generalized patterns $11$ and $1\mn1$. However, avoiding $1\mn1$
simply means having no repeated letters, so
\[
|[k]^n(1\mn1)|=\binom{k}{n}n!=(k)_n,
\]
the $n$-th lower factorial of $k$, which is 0 when $n>k$. Avoiding
$11$ is the same as having no repeated adjacent letters, so (see
Theorem \ref{tones})
\[
|[k]^n(11)|=k(k-1)^{n-1},
\]
for all $n\ge 1$. The remaining symmetry classes is those of
patterns $12$ and $1\mn2$. A word avoiding $1\mn2$ is just a
non-increasing string, so
\[
|[k]^n(1\mn2)|=\binom{n+k-1}{n}.
\]
Avoiding $12$ means there is no adjacent pair of letters in
increasing order, in other words, a string avoiding 12 is a
non-increasing string so (see Theorem \ref{t12})
\[
|[k]^n(12)|=\binom{n+k-1}{n}=|[k]^n(1\mn2)|.
\]
Indeed, the strings avoiding 12 are exactly those avoiding 1-2.

\begin{theorem}\label{tones}
Let $\tau=11\dots 1\in [1]^l$ be a generalized pattern. Then
\[
F_{\tau}(x;k;q)=\frac{1+(1-q)x\sum_{j=0}^{l-2} (kx)^j
-(1-q)(k-1)\sum_{i=2}^{l-1}x^d\sum_{j=0}^{l-1-d}
(kx)^j}{1-(k-1+q)x-(k-1)(1-q)(1-x^{l-2})\frac{x^2}{1-x}}.
\]
\end{theorem}
\begin{proof}
Let ${\langle j \rangle}_d=jj\dots j$, a string of $d$ letters
$j$. Then, by definition,
\[
\begin{split}
F_{\tau}(n,k;q;{\langle j \rangle}_d)  &=\sum_{\sigma\in
[k]^{n-d}} q^{\#\tau(\sigma,{\langle j \rangle}_d)}
            =\sum_{i=1}^k \sum_{\sigma\in [k]^{n-d-1}}
q^{\#\tau(\sigma,i,{\langle j \rangle}_d)}=\\
            &=\sum_{i\ne j} \sum_{\sigma\in [k]^{n-d-1}}
q^{\#\tau(\sigma,{\langle j \rangle}_d)}+ \sum_{\sigma\in
[k]^{n-d-1}} q^{\#\tau(\sigma,i,{\langle j \rangle}_{d+1})}.
\end{split}
\]
If we sum over all $j=1,2,\dots,k$, then for all $d\le l-2$ we
have
\[
\sum_{j=1}^k F_{\tau}(n,k;q;{\langle j
\rangle}_d)=(k-1)F_{\tau}(n-d;k;q)+ \sum_{j=1}^k
F_{\tau}(n,k;q;{\langle j \rangle}_{d+1}),
\]
hence
\[
F_{\tau}(n,k;q)=(k-1)\sum_{d=1}^{l-2} F_{\tau}(n-d,k;q) +G(n,k;q),
\]
where $G(n,k;q)=\sum_{j=1}^k \sum_{\sigma\in [k]^{n-(l-1)}}
q^{\#\tau(\sigma,{\langle j \rangle}_{l-1})}$. Again, by the same
above argument it is easy to see
\[
G(n,k;q)=(k-1)F_\tau(n-(k-l),k;q)+qG(n-1,k;q).
\]
Let we consider $F_\tau(n,k;q)-qF_\tau(n-1,k;q)$ with use the
above two equations, hence
\[
F_\tau(n,k;q)=(k-1+q)F_\tau(n-1,k;q)+(k-1)(1-q)\sum_{d=2}^{l-1}
F_\tau(n-d,k;q),
\]
for all $n\ge l$. Besides, $F_\tau(n,k;q)=k^n$ for all $n\le l-1$,
hence in terms of generating function the theorem holds.
\end{proof}

\begin{example} \label{t11}
\rm Theorem \ref{tones} yields for all $n\ge 1$,
\[
F_{11}(n,k;q)=k(q+k-1)^{n-1}.
\]
and
\[
F_{111}(x;k;q)=\frac{1+x(1+x)(1-q)}{1-(k-1+q)x-(k-1)(1-q)x^2}.
\]
Letting $q=1$, we get $|[k]^n(11)|=k(k-1)^{n-1}$, and
\[
\sum_{n\ge 0}|[k]^n(\underbrace{11\dots1}_l)|x^n=
\frac{1+x+\dots+x^{l-1}}{1-(k-1)x-\dots-(k-1)x^{l-1}}.
\]
For a different approach to the last formula, see \cite{Odlyzko},
Example 6.4 ff., pp. 1102--1103, for an easily generalizable case
of $k=2$.
\end{example}

We also obtain from the above that the number of strings in
$[k]^n$ with exactly $j$ occurrences of the generalized pattern
$11$ is $\binom{n-1}{j}k(k-1)^j$.

\begin{theorem}\label{t12}
For all $n\ge 2$ and $k\ge 1$,
\[
F_{12}(x;k;q)=\frac{1}{1+\frac{1-(1+(1-q)x)^k}{1-q}}.
\]
\end{theorem}
\begin{proof}
By definitions,
\[
\begin{split}
F_{12}(n,k;q;j) &=\sum_{i=1}^k \sum_{w\in [k]^{n-2}}
q^{\#12(wij)}=\\
        &=q\sum_{i=1}^{j-1} \sum_{w\in [k]^{n-2}}
q^{\#12(wi)}+\sum_{i=j}^{k} \sum_{w\in [k]^{n-2}}
q^{\#12(wi)}=\\
        &=(q-1)\sum_{i=1}^{j-1} \sum_{w\in [k]^{n-2}}
q^{\#12(wi)}+\sum_{i=1}^{k} \sum_{w\in [k]^{n-2}}
q^{\#12(wi)}=\\
        &=F_{12}(n-1,k;q)+(q-1)\sum_{i=1}^{j-1} \sum_{w\in
[k]^{n-2}} q^{\#12(wi)}.
\end{split}
\]
Besides, $F_{12}(n,k;q;1)=F_{12}(n-1,k;q)$, hence it is easy to see by
induction on $j$ that
\[
F_{12}(n,k;q;j)=\sum_{i=0}^{j-1} {\binom{j-1}{i} (q-1)^i
F_{12}(n-1-i,k;q)}.
\]
Therefore, for all $n\ge 1$,
\[
\begin{split}
F_{12}(n,k;q)&=\sum_{j=1}^k\sum_{i=0}^{j-1} {\binom{j-1}{i} (q-1)^i
F_{12}(n-1-i,k;q)}=\\
&=\sum_{j=1}^k {\binom{k}{j}(q-1)^{j-1}F_{12}(n-j,k;q)}.
\end{split}
\]
Besides, $F_{12}(0,k;q)=1$, hence the theorem holds for the
generating function $F_{12}(x;k;q)$.
\end{proof}

%===========================================================================

\section{Three-letter generalized patterns}

This section is divided into three subsections corresponding to
the three cases of $3$-letter generalized patterns: classical
patterns, patterns with exactly one adjacent pair of letters, and
patterns without internal dashes (i.e. with three consecutive
letters).

\subsection{Classical patterns}
The symmetry class representatives are $1\mn2\mn3$, $1\mn3\mn2$,
$1\mn1\mn2$, $1\mn2\mn1$, $1\mn1\mn1$. It is known \cite{Burstein}
that
\[
|[k]^n(1\mn2\mn3)|=|[k]^n(1\mn3\mn2)|=2^{n-2(k-2)}\sum_{j=0}^{k-2}{a_{k-2,j}\binom{n+2j}{n}},
\]
where
\[
a_{k,j}=\sum_{m=j}^{k}{C_m D_{k-m}}, \quad D_t=\binom{2t}{t},
\quad C_t=\frac{1}{t+1}\binom{2t}{t}.
\]
Recently, it was shown in \cite{BM} that
\[
|[k]^n(1\mn1\mn1)|=\sum_{i=0}^{k} {\binom{k}{i}\binom{i}{n-i}\frac{n!}{2^{n-i}}}
=\sum_{i=0}^{k}{B(i,n-i)(k)_i},
\]
where $(k)_i$ is the falling factorial, and $B(r,s)=\displaystyle
\frac{(r+s)!}{2^s (r-s)!s!}$ is the Bessel number of the first
kind. In particular, $f_{111}(n,k)=0$ when $n>2k$. It was also
proved in \cite{BM} that
\[
|[k]^n(1\mn2\mn1)|=|[k]^n(1\mn1\mn2)|=\sum_{j=0}^{k}{\binom{n+k-j-1}{n}c(n,n-j)},
\]
where $c(n,j)$ is the signless Stirling number of the first kind.

\subsection{Generalized patterns with exactly one adjacent pair of letters}
Let $F_\tau(x;k)$ be the generating function for the sequence
$\{f_\tau(n,k)\}_{n\ge 0}$ where $f_\tau(n,k)=|[k]^n(\tau)|$; that
is,
\[
F_\tau(x;k)=\sum_{n\ge 0} f_\tau(n,k)x^n.
\]
In the current subsection, we find explicit formulas or recurrence
formulas for $F_\tau(x;k)$ where $\tau$ is a three-letter
generalized pattern with exactly one adjacent pair of letters. The
symmetry class representatives are $11\mn1$, $11\mn2$, $21\mn1$,
$21\mn2$, $12\mn3$, $21\mn3$ and $13\mn2$. We will now find
$F_\tau(x;k)$ for each of the above patterns.

\begin{theorem}\label{11-1}
For all $k\ge 1$,
\[
\begin{split}
F_{11\mn1}(x;k)&=
\begin{vmatrix}
        B_k(x)      &   -A_k(x)     &   0   &   \dots   &   0 \\
        B_{k-1}(x)  &   1       &   -A_{k-1}(x) &   \dots   &   0 \\
        B_{k-2}(x)  &   0       &   1   &   \dots   &   0 \\
        \vdots      &   \vdots      &   \vdots   &   \ddots  &   \vdots   \\
        B_2(x)      &   0       &   0   &   \ddots  &   -A_{2}(x) \\
        B_1(x)+A_1(x)     &   0       &   0   &   \dots   &   1
\end{vmatrix}
\\
&=(-1)^{k-1}\prod_{i=1}^{k}{A_i(x)}+\sum_{j=1}^{k}
{\left((-1)^{k-j}B_j(x)\prod_{i=j+1}^{k}{A_i(x)}\right)}
\end{split}
\]
where $A_j(x)=\frac{jx^2}{1-(j-1)x}$ and
$B_j(x)=\frac{1+x}{1-(j-1)x}$.
\end{theorem}
\begin{proof}
Let $a_{n,k}=|[k]^n(11\mn1)|$, and let $a_{n,k}(i_1,\dots,i_d)$ be
the number of all words $\sigma\in [k]^n(11\mn1)$ such that
$\sigma_j=i_j$ for all $j=1,2,\dots d$. Then it follows that
$a_{n,k}=\sum_{i=1}^k a_{n,k}(i)$. On the other hand,
$a_{n,k}(i)=\sum_{j=1}^k a_{n,k}(i,j)$,
$a_{n,k}(i,j)=a_{n-1,k}(j)$ for $i\ne j$, and
$a_{n,k}(i,i)=a_{n-2,k-1}$. Therefore,
\[
a_{n,k}=k(a_{n-1,k}+a_{n-2,k-1})-a_{n-1,k}
\]
for all $n\ge 2$ and $k\ge 1$. Besides, $a_{1,k}=k$ and
$a_{0,k}=\delta_{0,k}$, hence, for all $k\ge 1$,
\[
F_{11\mn1}(x;k)=\frac{1+x}{1-(k-1)x}+\frac{kx^2}{1-(k-1)x}F_{11\mn1}(x;k-1),
\]
and $F_{11\mn1}(x;0)=1$. The rest follows by induction on $k$.
\end{proof}

\begin{theorem} \label{t11x2}\label{t22x1}
For all $k\ge 1$,
\[
F_{11\mn2}(x;k)=\prod_{j=0}^{k-1}{ \frac{1-(j-1)x}{1-(j+x)x} }.
\]
\end{theorem}
\begin{proof}
Let $\sigma\in [k]^n(11\mn 2)$; then there are $f_{11\mn2}(n,k-1)$
such words where $\sigma_j\ne k$ for all $j$. Let
$f_{11\mn2}(n,k;j)$ be the number of words
$\sigma\in[k]^n(11\mn2)$ where $\sigma_j=k$ and $j$ is minimal,
then $\sigma=\sigma_1 k \sigma_2$ for some
$\sigma_1\in[k-1]^{j-1}(11)$ and $\sigma_1\in[k]^{n-j}(11\mn2)$.
Moreover, for all $\sigma_1$ and $\sigma_2$ as above, we have
$\sigma=\sigma_1 k \sigma_2\in [k]^n(11\mn 2)$, so
\[
\begin{split}
f_{11\mn2}(n,k)&=f_{11\mn2}(n,k)+\sum_{j=1}^{n}{f_{11\mn2}(n,k;j)}\\
f_{11\mn2}(n,k;j)&=f_{11}(j-1,k-1)f_{11\mn2}(n-j,k)
\end{split}
\]
so
\[
f_{11\mn2}(n-1,k)+\sum_{j=2}^{n}{f_{11}(j-1,k-1)f_{11\mn2}(n-j,k)}.
\]
Hence, from the Example \ref{t11}, we get that, for all $n\ge 1$,
\[
f_{11\mn2}(n,k)=f_{11\mn2}(n,k-1)+f_{11\mn2}(n-1,k)+\sum_{j=0}^{n-2}
(k-1)(k-2)^j f_{11\mn2}(n-2-j,k).
\]
In addition, $f_{11\mn 2}(0,k)=\delta_{0,k}$, hence
\[
F_{11\mn2}(x;k)=\frac{1-(k-2)x}{1-(k-1)x-x^2}F_{11\mn2}(x;k-1).
\]
The rest follows easily.
\end{proof}

\begin{example} \label{112r0}
\rm For $k=1,2$, Theorem \ref{t11x2} yields
$F_{11\mn2}(x;1)=\frac{1}{1-x}$ and
$F_{11\mn2}(x;2)=\frac{1}{(1-x)(1-x-x^2)}$. In other words,
$f_{11\mn2}(n,1)=1$ and $f_{11\mn2}(n,2)=F_{n+3}-1$, where
$F_{n+3}$ is the $(n+3)$-rd Fibonacci number.
\end{example}

\begin{theorem}
Patterns $21\mn1$ and and $21\mn2$ are Wilf-equivalent, and we
have, for $k\ge 0$,
\[
F_{21\mn2}(x;k)=F_{21\mn1}(x;k)=1+\sum_{d=0}^{k-1} {\left( x^{d+1}
F_{21\mn1}(x;k-d)
\sum_{i=d}^{k-1}{(1-x)^{i-1}\binom{i}{d}}\right)}.
\]
\end{theorem}

\begin{proof}
Let us find $F_{21\mn1}(x;k)$ and $F_{21\mn2}(x;k)$ to prove that
the two are equal.

Let us derive the formula for $F_{21\mn1}(x;k)$.

Let $g_{n;k}(i_1,\dots,i_d)$ be the number of words $w$ in
$[k]^n(21\mn1)$ such that $w_j=i_j$ for all $j=1,2,\dots,d$ (in
other words, $w$ begins with the string $i_1\dots i_d$).

Consider $g_{n,k}(i,j)$. There are two cases. If $j\ge i$, then
the first $i$ cannot be part of any $21\mn1$ in a word $w$
beginning with $(i,j)$ and places no further restrictions on the
rest of the $w$, so $g_{n,k}(i,j)=g_{n-1,k}(j)$ if $j\ge i$. If
$j<i$, then deleting the first $i$ from $w$ starting with $(i,j)$,
we get a word $w'\in[k]^{n-1}(21\mn1)$ which contains exactly one
$j$, namely, as the first letter. Now let $w''\in[k-1]^{n-1}$ be
the word obtained from $w'$ by subtracting 1 from each letter
$>j$. Obviously, this is mapping is a bijection onto the set of
words in $[k-1]^{n-1}(21\mn1)$ starting with $j$ (since neither it
nor its inverse, i.e. adding one to each letter $\ge j$ except the
first letter, creates any new occurrences of $21\mn1$). Thus,
$g_{n,k}(i,j)=g_{n-1,k-1}(j)$ if $j<i$.

\begin{equation} \label{e21x1}
\begin{split}
g_{n,k}(i)&=\sum_{j=1}^{k}{g_{n,k}(i,j)}=\sum_{j=1}^{i-1}{g_{n-1,k-1}(j)}+\sum_{j=i}^{k}{g_{n-1,k}(j)}\\
&=f_{21\mn1}(n-1,k)+\sum_{j=1}^{i-1}{(g_{n-1,k-1}(j)-g_{n-1,k}(j))}.
\end{split}
\end{equation}

We can show by induction on $d$ that (\ref{e21x1}) implies
\[
g_{n,k}(d)=\sum_{j=0}^{d-1}{\binom{d-1}{j}\sum_{i=0}^{d-1-j}
{\binom{d-1-j}{i}(-1)^i f_{21\mn1}(n-1-i-j,k-j)}}.
\]
Hence for all $n\ge 1$,
\[
f_{21\mn1}(n,k)=\sum_{d=1}^{k}{\sum_{j=0}^{d-1}{\binom{d-1}{j}\sum_{i=0}^{d-1-j}{
\binom{d-1-j}{i}(-1)^i f_{21\mn1}(n-1-i-j,k-j)}}}
\]
with $f_{21\mn1}(0,k)=1$.

If we multiply through by $x^n$ and sum over all $n\ge 1$, we get
\[
F_{21\mn1}(x;k)-1=\sum_{d=1}^{k}
{\sum_{j=0}^{d-1}{\binom{d-1}{j}x^{j+1}(1-x)^{d-1-j}
F_{21\mn1}(x;k-j)}},
\]
which means that
\[
F_{21\mn1}(x;k)-1=\sum_{d=0}^{k-1}{\left(x^{d+1}F_{21\mn1}(x;k-d)
\sum_{i=d}^{k-1}{(1-x)^{i-1}\binom{i}{d}}\right)}.
\]

Now let us find the formula for $F_{21\mn2}(x;k)$.

Similarly, we define $h_{n,k}(i_1,\dots,i_d)$ to be the number of
words $w$ in $[k]^n(21\mn2)$ such that $w_j=i_j$ for all
$j=1,2,\dots,d$. As before, $h_{n,k}(i,j)=h_{n-1,k}(j)$ if $j\ge
i$. If $j<i$ and $w\in[k]^n(21\mn2)$ starts with $(i,j)$, then $w$
must contain only one $i$. Thus, deleting $i$ and subtracting 1
from each letter $>i$ is a bijection from the set of words in
$[k]^n(21\mn2)$ that start with $(i,j)$ onto the set of words in
$[k-1]^{n-1}(21\mn2)$ that start with $j$. Therefore,
$h_{n,k}(i,j)=h_{n-1,k-1}(j)$ if $j<i$, so it is easy to see that
$h_{n,k}(i)$ satisfies the same recurrence \ref{e21x1} as
$g_{n,k}(i)$, i.e.
\[
h_{n,k}(i)=f_{21\mn2}(n-1,k)+\sum_{j=1}^{i-1}{(h_{n-1,k-1}(j)-h_{n-1,k}(j))}.
\]
Noting that $g_{n,k}(i)$ and $h_{n,k}(i)$ are the same when
$n=0,1$ or $k=0,1$, we see that $F_{21\mn1}(x;k)=F_{21\mn2}(x;k)$
for $k\ge 0$.
\end{proof}

\begin{example}
\[
\begin{split}
F_{21\mn2}(x;0)&=F_{21\mn1}(x;0)=1\\
F_{21\mn2}(x;1)&=F_{21\mn1}(x;1)=\frac{1}{1-x}\\
F_{21\mn2}(x;2)&=F_{21\mn1}(x;2)=\frac{1}{(1-x)^2}+\frac{x^2}{(1-x)^3}\\
F_{21\mn2}(x;3)&=F_{21\mn1}(x;3)=\frac{1-3x+6x^2-5x^3+3x^4-x^5}{(1-x)^6}
\end{split}
\]
\end{example}

%t12x3==t32x1==t1x23==t3x21==t21x3==t3x12==t23x1==t1x32

\begin{theorem}\label{t12x3}
Patterns $12\mn3$ and $21\mn3$ are Wilf-equivalent, and for all
$k\ge 1$,
\[
F_{12\mn3}(x;k)=F_{21\mn3}(x;k)=\prod_{j=0}^{k-1}
\frac{1}{1-\frac{x}{(1-x)^j}}.
\]
\end{theorem}
\begin{proof}
Let $\sigma=(\sigma',k,\sigma'')\in [k]^n$ such that $\sigma_j=k$
and $j$ is minimal. Hence $\sigma$ avoids $12\mn3$ (respectively,
$21\mn3$) if and only if $\sigma'\in [k-1]^{j-1}$ avoids $12$
(respectively, $21$) and $\sigma''\in [k]^{n-j}$ avoids $12\mn3$
(respectively, $21\mn3$). Theorem \ref{t12} yields for all $n\ge
0$ and $k\ge 2$,
\[
\begin{split}
f_{12\mn3}(n,k)&=f_{12\mn3}(n,k-1)+\sum_{j=1}^{n}{\binom{j-1+k-2}{k-2}}
f_{12\mn3}(n-j,k),\\
f_{21\mn3}(n,k)&=f_{21\mn3}(n,k-1)+\sum_{j=1}^{n}{\binom{j-1+k-2}{k-2}}
f_{21\mn3}(n-j,k).
\end{split}
\]
Since $f_{12\mn3}(n,1)=f_{21\mn3}(n,1)=1$ as well, we get that
$F_{12\mn3}(x;k)=F_{21\mn3}(x;k)$ and
\[
F_{12\mn3}(x;k)=F_{12\mn3}(x;k-1)+\frac{x}{(1-x)^{k-1}}F_{12\mn3}(x;k).
\]
\end{proof}

%=t13x2==========
\begin{lemma}\label{e13x2}
Let $a_{n,k}(i)$ be the number of words $\sigma\in [k]^n(13\mn2)$
such that $\sigma_1=i$. Then, for $n\ge 1$, $k\ge 0$
\[
f_{13\mn2}(n,k)=\sum_{j=1}^k a_{n,k}(j),\quad
a_{n,k}(j)=\sum_{i=1}^{j+1}{a_{n-1,k}(i)}+\sum_{i=j+1}^{k-1}
{a_{n-1,i}(j+1)}.
\]
In addition, $f_{13\mn2}(n,1)=1$ for all $n$.
\end{lemma}
\begin{proof}
The first identity holds by definition, and
\[
a_{n,k}(j)=\sum_{i=1}^{k}{a_{n,k}(j,i)}=\sum_{i=1}^{j+1}{a_{n,k}(j,i)}
+\sum_{i=j+2}^{k}{a_{n,k}(j,i)},
\]
where $a_{n,k}(j,i)$ is the number of $\sigma\in [k]^n(13\mn 2)$
such that $\sigma_1=j$ and $\sigma_2=i$. If $i\le j+1$, then the
first letter $j$ places no restriction on the rest of $\sigma$, so
$a_{n,k}(j,i)=a_{n-1,k}(i)$. If $j+2\le i\le k$, then the rest of
$\sigma$ (positions 3 to $n$) may not contain any letters from
$j+1$ to $i-1$ (for a total of $i-j-1$ letters), so
$a_{n,k}(j,i)=a_{n-1,k+j+1-i}(j+1)$. Hence
\[
a_{n,k}(j)=\sum_{i=1}^{j+1}{a_{n-1,k}(i)}+\sum_{i=j+2}^{k}
{a_{n-1,k+j+1-i}(j+1)}.
\]
The rest is easy to obtain.
\end{proof}

We also note that $a_{n,k}(k-1)=a_{n,k}(k)=f(n-1,k)$ and that the
second formula in Lemma \ref{e13x2} implies
\begin{multline*}
a_{n,k}(j)-a_{n,k}(j-1)-a_{n,k-1}(j)+a_{n,k-1}(j-1)=\\=
\begin{cases}
a_{n-1,k}(j+1)-a_{n-1,k-1}(j), \quad \text{for $2\le j\le k-1$},\\
a_{n-1,k}(j+1)-a_{n-1,k-1}(j)+a_{n-1,k}(j), \quad \text{for $j=1$},\\
0, \quad \text{for $j=k$}.
\end{cases}
\end{multline*}

From Lemma \ref{e13x2}, it is easy to obtain explicit formulas for
small values of $k$.

\begin{theorem}
\[
\begin{split}
F_{13\mn2}(x;1)&=\frac{1}{1-x},\\
F_{13\mn2}(x;2)&=\frac{1}{1-2x},\\
F_{13\mn2}(x;3)&=\frac{(1-x)^2}{(1-2x)(1-3x+x^2)},\\
F_{13\mn2}(x;4)&=\frac{1-4x+6x^2-3x^3}{(1-3x)(1-2x)(1-3x+x^2)}.
\end{split}
\]
\end{theorem}

%=======================

\subsection{Generalized Patterns without internal dashes}
The symmetry class representatives are $111$, $122$, $212$, $123$,
$213$. In the current subsection, we find explicit formulas for
$F_\tau(x;k)$ for each of these representatives $\tau$. Example
\ref{t11} yields the following result for $111$.

\begin{theorem}
For all $k\ge 1$,
\[
F_{111}(x;k)=\frac{1+x+x^2}{1-(k-1)x-(k-1)x^2}.
\]
\end{theorem}

\begin{theorem}\label{t122}
\[
F_{122}(x;k)=\frac{x}{(1-x^2)^k-(1-x)}, \quad k\ge 0.
\]
\end{theorem}

\begin{proof}
Let $\sigma\in [k]^n(122)$ where $\sigma$ contains $j$ letters
$k$. If $j=0$, then there are $f_{122}(n,k-1)$ such words. Let
$j>1$, $\sigma=(\sigma',k,\sigma'')$, and let $\sigma_{r}=k$ where
$r$ minimal (so $\sigma'$ does not contain $k$). If
$\sigma_{r+1}\ne k$ then $\sigma$ avoids $122$ if and only if
$\sigma'\in [k-1]^{r-1}(122)$, and $\sigma''\in [k]^{n-r}(122)$
such that $\sigma''_1\ne k$, so there are
\[
\sum_{r=1}^{n}f_{122}(r-1,k-1)(f_{122}(n-r,k)-f_{122}(n-1-r,k))
\]
such words. If $\sigma_{r+1}=k$, then $\sigma'=\emptyset$, and
$\sigma$ avoids $122$ if and only if $\sigma''$ avoids $122$, so
there are $f_{122}(n-2,k)$ such words. Hence
\begin{multline*}
f_{122}(n,k)=f_{122}(n,k-1)+f_{122}(n-2,k)\\
+\sum_{i=1}^{n}{f_{122}(i-1,k-1)(f_{122}(n-i,k)-f_{122}(n-1-i,k))}
\end{multline*}
for all $n\ge 2$ and $k\ge 1$, therefore,
\begin{multline*}
f_{122}(n,k)=f_{122}(n,k-1)+f_{122}(n-2,k)\\
+\sum_{i=0}^{n-1}{f_{122}(i,k-1)f_{122}(n-1-i,k)}
-\sum_{i=0}^{n-1}{f_{122}(i,k-1)f_{122}(n-2-i,k))}.
\end{multline*}
In addition, $f_{122}(1,1)=f_{122}(0,1)=f_{122}(0,0)=1$ and
$f_{122}(n,k)=0$ for $n<0$, hence, taking the generating
functions, we see that
\[
F_{112}(x;k)=F_{112}(x;k-1)+x^2F_{112}(x;k)+(x-x^2)F_{112}(x;k)F_{112}(x;k-1).
\]
Now, after dividing through by $F_{112}(x;k)F_{112}(x;k-1)$ and
some routine manipulations, we see that the theorem holds.
\end{proof}

\begin{example}
\rm Theorem \ref{t122} yields
$F_{122}(x;2)=\frac{1}{(1-x)(1-x-x^2)}$, which means
$f_{122}(n,2)=F_{n+3}-1$, where $F_{n+3}$ is the $(n+3)$-rd
Fibonacci number.
\end{example}

\begin{theorem}
\[
F_{212}(x;k)=\frac{1}{1-x\sum_{j=0}^{k-1}{\frac{1}{1+jx^2}}},
\quad k\ge 0.
\]
\end{theorem}
\begin{proof}
Let $d_{212}(n,k)$ be the number of words $\sigma\in
[k]^{n+1}(212)$ such that $\sigma_1=k$, and let $D_{212}(x;k)$ be
the generating function for $d_{212}(n,k)$, that is
$D_{212}(x;k)=\sum_{n\ge 0}{d_{212}(n,k)x^n}$.

Let $\sigma=(\sigma',k,\sigma'')\in [k]^n(212)$ so that $\sigma'$
does not contain $k$. If $\sigma=\sigma'$, then there are
$f_{212}(n,k-1)$ such words. Otherwise, $\sigma'\in [k-1]^{j-1}$
for some $j=1,\dots,n$, so there are
$f_{212}(j-1,k-1)d_{212}(n-j,k)$ such words. Therefore, for all
$n\ge 1$,
\[
f_{212}(n,k)=f_{212}(n,k-1)+\sum_{j=1}^{n}{f_{212}(j-1,k-1)d_{212}(n-j,k)}.
\]
In addition, $f_{212}(0,k)=0$ for all $k\ge 1$, hence
\[
F_{212}(x;k)=(1+xD_{212}(x;k))F_{212}(x;k-1).
\]

Now let $\sigma=(k,\sigma'')\in [k]^{n+1}(212)$. The first $k$ is
not part of any occurrence of $212$, so we can delete it to get
any word $\sigma''\in [k]^n(212)$. We can obtain an occurrence of
$212$ by adding the first $k$ back only if
$\sigma=(k,i,k,\sigma''')$ for some $i<k$ and $(k,\sigma''')\in
[k]^{n-1}(212)$. Then neither the second nor the third letter of
$\sigma$ can start $212$, so the number of such ``bad'' words is
$(k-1)d_{212}(n-2,k)$. Hence,
\[
d_{212}(n,k)=f_{212}(n,k)-(k-1)d_{212}(n-2,k), \quad n\ge 2.
\]
Taking generating functions of both sides, we obtain
\[
D_{212}(x;k)=F_{212}(x;k)-(k-1)x^2D_{212}(x;k).
\]
Now, solving the two generating function recurrences above, we get
\[
\frac{1}{F_{212}(x;k)}=\frac{1}{F_{212}(x;k-1)}-\frac{x}{1+(k-1)x^2},
\]
which implies the theorem.
\end{proof}

\begin{theorem}\label{t123}
\[
F_{123}(x;k)=\frac{1}{\sum_{j=0}^{k}{a_j\binom{k}{j}x^j}}, \quad
k\ge 0,
\]
where $a_{3m}=1$, $a_{3m+1}=-1$, $a_{3m+2}=0$ for all $m\ge 0$.
\end{theorem}
\begin{proof}
By definition, $f_{123}(n;k)=k^n$ for $k=0,1,2$, so the formulas
for $F_{123}(x;k)$ hold for $k=0,1,2$. Let $d_{123}(n,k)$ be the
number of words $\sigma\in [k]^n(123)$ such that $(\sigma,k+1)$
also avoids 123, and let $D_{123}(x;k)$ be the generating function
for $d_{123}(n,k)$ with $k$ fixed.

\begin{lemma}\label{l1}
\[
D_{123}(x;3)=\frac{1-3x^2+x^3}{1-3x+x^3},\quad
F_{123}(x;3)=\frac{1}{1-3x+x^3}.
\]
\end{lemma}
\begin{proof}
Let $\sigma\in [3]^n(123)$, and let $f_{123}(n,k;a_1,\dots,a_d)$
be the number of words $\sigma\in [3]^n(123)$ such that
$\sigma_i=a_i$ for all $i=1,2\dots,d$. Deleting the first letter
$i\in\{1,2,3\}$ of $\sigma$, we get a word
$\sigma'\in[3]^{n-1}(123)$. Adjoining $i$ in front of
$\sigma'\in[3]^{n-1}(123)$, we get a word $\sigma\in[k]^n(123)$ or
$\sigma=(1,2,3,\sigma'')$ for some $\sigma''\in[3]^{n-3}(123)$.
Hence,
\[
f_{123}(n,3)=3f_{123}(n-1,3)-f_{123}(n-3,3), \quad n\ge 3.
\]
Besides, $f_{123}(n,3)=3^n$ for $n=0,1,2$, hence the formula for
$F_{123}(x;3)$ holds.

Similarly, we have
\[
d_{123}(n,3)=3d_{123}(n-1,3)-d_{123}(n-3,3), \quad n\ge 4.
\]
Besides, $d_{123}(0,3)=1$, $d_{123}(1,3)=3$, $d_{123}(2,3)=6$, and
$d_{123}(3,3)=18$, hence the formula $D_{123}(x;3)$ holds as well.
\end{proof}

Now, we are ready to prove the two main recurrences.

\begin{lemma}\label{l2}
For all $k\ge 4$,
\[
F_{123}(x;k)=\frac{F_{123}(x;k-1)}{1-xD_{123}(x;k-1)},\quad
D_{123}(x;k)=\frac{x+(1-x)D_{123}(x;k-1)}{1-xD_{123}(x;k-1)}.
\]
\end{lemma}
\begin{proof}
Let $\sigma=(\sigma',k,\sigma'')\in [k]^n(123)$ be such that
$\sigma_j=k$ and $\sigma'\in [k-1]^{j-1}$. If $\sigma=\sigma'$,
then there are $f_{123}(n,k-1)$ such words, otherwise, there are
$d_{123}(j-1,k-1)f_{123}(n-j,k)$ of them. Therefore,
\[
f_{123}(n,k)=f_{123}(n,k-1)+
\sum_{j=1}^{n}{d_{123}(j-1,k-1)f_{123}(n-j,k)},
\]
for all $n\ge 1$. Besides, $f_{123}(0,k)=0$ for $k\ge 1$, hence,
taking generating functions, we get that
\[
F_{123}(x;k)=\frac{F_{123}(x;k-1)}{1-xD_{123}(x;k-1)}.
\]
Similarly, we have for all $n\ge 1$,
\[
d_{123}(n,k)=d_{123}(n,k-1)+\sum_{i=1}^{n-1}{d_{123}(i-1,k-1)d_{123}(n-i,k)},
\]
or, equivalently,
\[
d_{123}(n,k)=d_{123}(n,k-1)-d_{123}(n-1,k-1)+\sum_{i=1}^{n}{d_{123}(i-1,k-1)d_{123}(n-i,k)}
\]
for all $n\ge 1$. Besides, $d_{123}(0,k)=1$, and $d_{123}(n,k)=0$
for $n<0$, hence, taking generating functions, we get
\[
D_{123}(x;k)=\frac{x+(1-x)D_{123}(x;k-1)}{1-xD_{123}(x;k-1)}.
\]
\end{proof}

Finally, Lemmas \ref{l1} and \ref{l2} together yield us Theorem
\ref{t123}. We also note that the same Lemmas yield that
\[
D_{123}(x;k)=\frac{\sum_{j=0}^{k}{b_j\binom{k}{j}x^j}}{\sum_{j=0}^{k}{a_j\binom{k}{j}x^j}},
\quad k\ge 0,
\]
where $b_{3m}=a_{3m}=1$, $b_{3m+1}=a_{3m+2}=0$,
$b_{3m+2}=a_{3m+1}=-1$ for all $m\ge 0$.
\end{proof}

\begin{theorem}
\[
F_{213}(x;k)=\frac{1}{1-x-x\sum_{i=0}^{k-2}{\prod_{j=0}^{i}{(1-jx^2)}}},
\quad k\ge 1,
\]
and $F_{213}(x;0)=1$.
\end{theorem}
\begin{proof}
For $k=0$ the theorem is trivial, so we may assume $k\ge 1$. Let
$d_{213}(n,k)$ be the number of all $\sigma\in [k]^n(213)$ such
that $(\sigma,k+1)$ also avoids 213, and let $D_{213}(x;k)$ be the
generating function for $d_{213}(n,k)$ with $k$ fixed.

Similarly to Lemma \ref{l2} we have
\[
f_{213}(n,k)=f_{213}(n,k-1)+\sum_{j=1}^{n}{d_{213}(j-1,k-1)f_{213}(n-j,k)},
\quad n\ge 1.
\]
Besides $f_{213}(0,k)=0$ for $k\ge 1$, hence, taking generating
functions, we get that
\[
F_{213}(x;k)=\frac{F_{213}(x;k-1)}{1-xD_{213}(x;k-1)}.
\]
Similarly to Lemma \ref{l2}, we have for all $n\ge 1$,
\begin{multline*}
d_{213}(n,k)=d_{213}(n,k-1)+d_{213}(n-1,k-1)+d_{213}(n-2,k-1)\\
+\sum_{i=1}^{n-2}{d_{213}(i-1,k-1)d_{213}(n-i,k)},
\end{multline*}
or, equivalently,
\begin{multline*}
d_{213}(n,k)=d_{213}(n,k-1)-(k-1)d_{213}(n-2,k-1)\\
+\sum_{i=1}^{n}{d_{213}(i-1,k-1)d_{213}(n-i,k)},
\end{multline*}
for all $n\ge 1$. Besides, $d_{213}(0,k)=1$, and $d_{213}(n,k)=0$
for $n<0$, hence, taking generating functions, we obtain
\[
D_{213}(x;k)=\frac{(1-(k-1)x^2)D_{213}(x;k-1)}{1-xD_{213}(x;k-1)}.
\]
\end{proof}

%===========================================================================

\end{document}